\documentclass{article}

\usepackage{epsfig}
\usepackage{subfigure}
\usepackage{calc}
\usepackage{graphicx}

\usepackage{amsfonts}
\usepackage{amsmath}
\usepackage{amssymb}
\usepackage{latexsym}

\newcommand{\adots}{\mathinner{\mkern2mu\raise1pt\hbox{.}\mkern2mu%
\raise4pt\hbox{.}\mkern2mu\raise7pt\hbox{.}\mkern1mu}}

\begin{document}

\title{A Joint Criterion for Reachability and
Observability of Nonuniformly Sampled Discrete
Systems}
\date{}
\author{Amparo F\'{u}ster-Sabater\\
{\small (1) Instituto de Electr\'{o}nica de Comunicaciones, C.S.I.C.}\\
{\small Serrano 144, 28006 Madrid, Spain} \\
{\small amparo@iec.csic.es}}

\maketitle

\begin{abstract}

A joint characterization of reachability (controllability) and
observability (constructibility) for linear SISO nonuniformly
sampled discrete systems is presented. The work generalizes to the
nonuniform sampling the criterion known for the uniform sampling.
Emphasis is on the nonuniform sampling sequence, which is believed
to be an additional element for analysis and handling of discrete
systems.

\end{abstract}

\section{Introduction}
\footnotetext{Work supported by Ministerio de Educaci\'{o}n y
Ciencia (Spain), Project TIC-0386.\\
IEEE Transactions on Automatic Control. Volume 36, No. 11, pp. 1281-1284. Nov. 1991.\\
DOI:10.1109/9.100938} The concepts of controllability and observability, first introduced
by Kalman \cite{Kalman1}, are still two fundamental questions in modern
control theory. On the other hand, the enormous increase in the use
of digital computers has stimulated studies in the field of discrete
systems. In particular, the problem of controllability (reachability)
and observability of discrete systems has been already treated in the
literature in a generalized form \cite{Ackermann}, \cite{Astrom}, \cite{Kailath} - \cite{Troch}.

Most of the previous references are only concerned with uniformly
sampled discrete systems. However, the general case of nonuniform
sampling offers a wider range of situations in the analysis of
these concepts for discrete systems. The present note tackles this
problem and, in fact, a joint characterization of reachability and
observability for linear single-input single-output nonuniformly
sampled discrete systems has been developed. Under several
conditions, a right choice of the sampling instants would
guarantee the aforementioned internal properties.

This note emphasizes the importance of the nonuniform sampling
sequence against other system parameters. Besides the above
considerations, nonuniform sampling is believed to be an
alternative solution to relevant problems such as propagation of
measuring and/or rounding errors, communication delays in complex
computer-controlled systems, identification, etc.

\section{A Joint Criterion for Reachability and Observability of Nonuniformly Sampled Discrete Systems}

\subsection{Model Description}
Consider a linear time-invariant SISO dynamic system
\begin{equation*}
\dot{X}(t) = AX(t) + bu(t) \qquad X_0 = X(0)
\end{equation*}
\begin{equation}\label{eq:1}
Y(t) = c X(t)
\end{equation}

where $X \in \mathbb{R}^n$ denotes the state vector and $u, y \in \mathbb{R}$ are the scalar
input and output, respectively. $A(n \times n)$, $b(n \times 1)$, $c(1 \times n)$ are
real and constant matrices and $n$ is the order of the system.

As additional assumption, the realization $( A , b, c)$ is
required to be minimal.

\subsection{The Reachability Problem}
Let $(A , b, c)$ be an arbitrary minimal realization of order $n$ for the kind of system under study.
The solution of the state-space equation at the sampling instant $t_n$ can be written as

\begin{equation}\label{eq:2}
X(t_n)= exp(A\,t_n)X_0 + [G_{n-1}, \ldots , G_0]
\left[
  \begin{array}{c}
    u_{n-1} \\
    $\vdots$ \\
    u_0 \\
  \end{array}
\right]
\end{equation}

where
\begin{equation}\label{eq:3}
G_i = exp(A(t_n - t_i))b \qquad (i=0, \ldots, n-1)
\end{equation}
and $u_j$ is the value of the impulse input at time $t_i$.

In mathematical terms, the condition of \textit{n}-reachability
will be accomplished if and only if the matrix $[G_{n-1},  \ldots
, G_0]$ has full rank.

Vectors $G_i$ can be rewritten as
\begin{equation}\label{eq:4}
G_i = B \, exp(J(t_n - t_i))y_0
\end{equation}
where $J$ is the Jordan canonical form of the matrix $A$, $B$ is the
invertible matrix of the change of basis, and

\begin{equation}\label{eq:5}
y_0 = B^{-1}b.
\end{equation}

Therefore, in order to guarantee the \textit{n}-reachability we
must compute the value of

\begin{equation}\label{eq:6}
det[\,exp(J \alpha_m)y_0\,] \qquad (m =0, \ldots, n-1)
\end{equation}
with

\begin{equation}\label{eq:7}
\alpha_m = t_{n-1} - t_{n-m-1} \qquad (\alpha_0 = 0)
\end{equation}
and ensure that such a determinant is nonnull.

First, we denote the components of $y_0$ by
\begin{equation}\label{eq:8}
y_0 = (y_{1}^1, \ldots, y_{m1}^1, \ldots, y_{1}^r, \ldots, y_{mr}^r)'
\end{equation}
where $m_j \;(j = 1, \ldots , r)$ is the multiplicity of the $r$ different eigenvalues
of the matrix $A$ with $r \leq n$. ($'$ denotes the transpose).

Then, we denote the reachability canonical form \cite{Kailath} of
this arbitrary minimal realization by $(A_{re}, b_{re}, c_{re})$.
Finally, the impulse response \cite{Kailath} for this kind of
system can be written as

\begin{equation}\label{eq:9}
h(t) = c \, exp(At) b=\sum\limits_{i=1}^{n} c_i \varphi_i(t)
\end{equation}

where $c_i \in \mathbb{C}$ are constant coefficients and
$\varphi_i; \mathbb{R} \rightarrow \mathbb{C} \;\; (i = 1,\ldots ,
n)$ is the fundamental system of solutions of an
\textit{n}th-order homogeneous linear differential equation.

Now, making use of the Laplace's expansion by minors, the
determinant of (\ref{eq:6}) can be factorized as follows:

\begin{equation}\label{eq:10}
det[\,exp(J \alpha_m)y_0\,]= N_1 N_2 \, det[\,\varphi_i(\alpha_m)\,].
\end{equation}

Let us consider separately each one of these factors
\begin{equation}\label{eq:11}
N_1= \frac{1}{0!}\ldots \frac{1}{(m_1-1)!} \ldots \frac{1}{0!}\ldots \frac{1}{(m_r-1)!}
\end{equation}

The term $N_1$ is related to the multiplicity of the eigenvalues
of the matrix $A$ and will always be nonnull

\begin{equation}\label{eq:12}
N_2=det \left[%
\begin{array}{ccc}
  \left[%
\begin{array}{ccc}
  y_1^1 & \ldots & y_{m1}^1 \\
  \vdots & \adots & \, \\
  y_{m1}^1 & \, & \, \\
\end{array}%
\right] & \, & \, \\
  \, & \ddots & \, \\
  \, & \, & \left[%
\begin{array}{ccc}
  y_1^r & \ldots & y_{mr}^r \\
  \vdots & \adots & \, \\
  y_{mr}^r & \, & \, \\
\end{array}%
\right] \\
\end{array}%
\right]
\end{equation}

By similarity transformations
\begin{equation}\label{eq:13}
y_0 = B^{-1} b = B_{re}^{-1} b_{re}= (C_1,C_2, \ldots, C_n)'
\end{equation}

where $B_{re}$ is the matrix of the change the basis of $A_{re}$
to the Jordan canonical form. The term $N_2$ is related to the
weighting coefficients $(C_i)$ of the characteristic modes in
(\ref{eq:9}). Remark that $N_2$ will be nonnull if and only if
\begin{equation}\label{eq:14}
y_{mj}^j \neq 0 \qquad (j = 1, \ldots, r).
\end{equation}

This holds, according to the previous meaning of the components of
$y_0$, because only minimal realizations are considered.

Finally, $[\varphi_i(\alpha_m)] \;\; (i= 1, \ldots , n ; \,m = 0,
\ldots , n - 1)$ is a $n \times n$ matrix involving characteristic
modes and sampling instants. From (\ref{eq:10}), the following
result is derived.

\textit{Lemma 1:} An arbitrary minimal realization $(A, b, c)$ is
completely \textit{n}-reachable (reachable in $n$ steps) if and only if $n$ consecutive
sampling instants are chosen in such a way that
\begin{equation}\label{eq:15}
det[\,\varphi_i(\alpha_m)\,]\neq 0 \qquad (i = 1, \ldots, n;\, m = 0, \ldots, n-1).
\end{equation}
\textit{Proof:} The proof is evident from the factorization of
(\ref{eq:10}).

Note that if the scalar input had been defined as a control of the
form
\begin{equation}\label{eq:16}
u(t) = u(t_i) = u_i \qquad t_i \leq t < t_{i+1}
\end{equation}

which is more realistic, then the above result would still be
valid. In fact, the presence of a data hold does not affect the
system characteristic modes.

\subsection{General Criterion}
By duality, the \textit{n}-observability characterization is straightforwardly
derived. Indeed, a similar lemma (concerning the \textit{n}-observability)
can be stated by just changing the term reachability to
observability in Lemma 1.

The results obtained previously can be unified as follows.

\textit{Theorem:} An arbitrary minimal realization $(A, b, c)$ for the
kind of system under study is jointly \textit{n}-reachable and \textit{n}-observable
if and only if $n$ consecutive sampling instants, not necessarily
equidistant, are chosen in such a way that
\begin{equation}\label{eq:17}
det[\,\varphi_i(\alpha_m)\,]\neq 0 \qquad (i = 1, \ldots, n; \, m = 0, \ldots, n-1).
\end{equation}
\textit{Proof:} The proof is evident from the two previous lemmas.
It must be noticed that the condition (17) depends exclusively on
the system characteristic modes and the sampling instants.

Condition (17) imposes a rather weak restriction for the choice of
the sampling instants. In fact, time intervals can be specified so that
complete reachability and observability are preserved. Particular to
the uniform case, (17) becomes the condition imposed on the
sampling interval $T$ (see \cite{Kailath}, \cite{Kalman2}) when a uniform sampling is
considered.

We can also remark that, for this kind of system, the \textit{n}-reachability
and \textit{n}-observability are inseparable concepts. These systems are
either \textit{n}-reachable and \textit{n}-observable or, if not, they are neither
\textit{n}-reachable nor \textit{n}-observable.

On the other hand, the \textit{n}-controllability
(\textit{n}-constructibility) of an arbitrary minimal realization
can be considered as a weakening of the condition of
\textit{n}-reachability (\textit{n}-observability) \cite{Kailath}.
Therefore, the pair
\textit{n}-controllability/\textit{n}-constructibility can be
characterized as a corollary of the preceding theorem. Indeed,
both properties will be accomplished if and only if

\begin{equation}\label{eq:18}
exp(At_n)\,X_0 \in \mathbb{R} [G_{n-1}, \ldots , G_0]
\end{equation}

where $\mathbb{R}[ \ldots ]$ denotes the range space of the columns $G_i$.

From the previous theorem and following a similar reasoning, the
next corollary can be proved.

\textit{Corollary:} An arbitrary minimal realization $(A, b, c)$ for the
kind of system under study is jointly \textit{n}-controllable and \textit{n}-constructible
if and only if $n + 1$ sampling instants, not necessarily equidistant,
are chosen in such a way that
\begin{equation}\label{eq:19}
(\varphi_1(\alpha_n), \ldots, \varphi_n(\alpha_n))' \in \mathbb{R}[(\varphi_1(\alpha_n), \ldots, \varphi_n(\alpha_n))']
\end{equation}

(with $\alpha_m$, $\varphi_i$ defined as before) $\;(m = 0, \ldots
, n - 1)$
\begin{equation}\label{eq:20}
\alpha_n = t_n - t_0.
\end{equation}

Remark that the condition of \textit{n}-controllability /
\textit{n}-constructibility involves one sampling instant more
than in the preceding characterization. Note also that the pair
\textit{n}-reachability / \textit{n}-observability implies the
pair \textit{n}-controllability/ \textit{n}-constructibility, but
the converse may not be true.

\begin{figure}
\begin{center}
\includegraphics[bb= 10 12 330 210, scale=1.2]{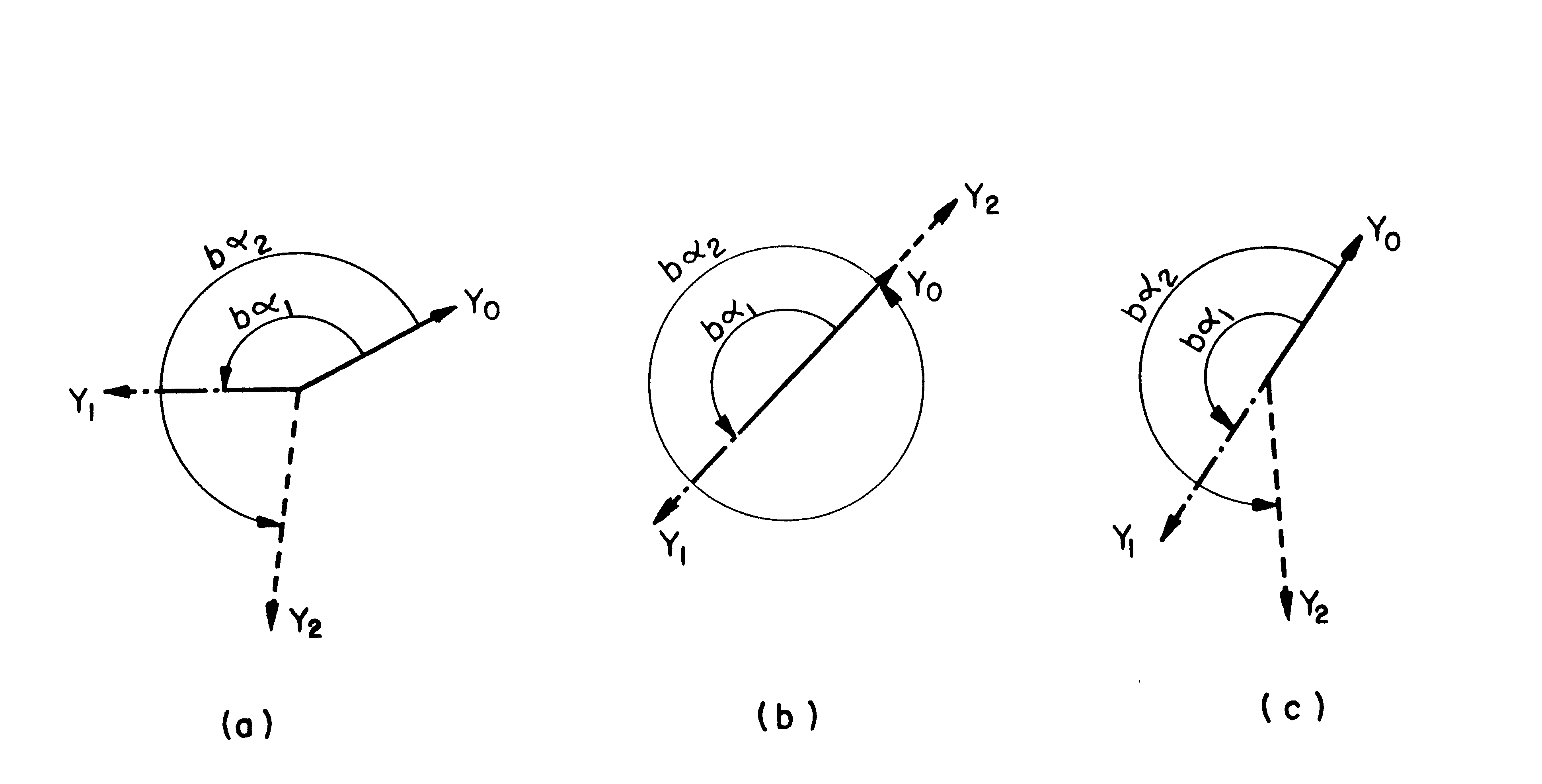}
\caption{Relationships among the vectors $Y_0, Y_1, Y_2$}
\label{figure:headings1}
\end{center}
\end{figure}

\section{A Simple Strategy for Choosing Sampling Instants}
In order to give more insight on the restriction imposed by (17), a
simple example is presented. Consider a 2nd-order realization where
the matrix $A$ has a pair of complex eigenvalues
\begin{equation}\label{eq:21}
a + jb \qquad (b > 0).
\end{equation}
Now, (6) can be written as
\begin{equation}\label{eq:22}
det[\,Y_0, Y_1\,]= det[\,exp(J\alpha_0)y_0, \, exp(J\alpha_1)y_0\,]
\end{equation}
where
\begin{equation}\label{eq:23}
J= \left[%
\begin{array}{cc}
  a & -b \\
  b & a \\
\end{array}%
\right]
\end{equation}

and $\alpha_i, y_0$ are defined as before.

A geometric interpretation is very easy. Indeed, the generic
operator $exp(J\alpha)$ applied on the vector $y_0$ can be viewed as
follows:
\begin{enumerate}
\item a counter-clockwise rotation through $b\alpha$ rad;
\item a stretching (or shrinking) of the length of $y_0$ by a
    factor $exp(a\alpha)$
\end{enumerate}

Thus, a \textit{necessary and suficient condition} to preserve
complete reachability and observability is that

\begin{equation}\label{eq:24}
b(t_1 - t_0) \neq k\Pi \qquad k=0, 1, \ldots \;.
\end{equation}

Remark that given a first and arbitrary sampling instant $t_0$, just
point values on the time axis will be forbidden. Thus, (17) is
actually not very restrictive.

On the other hand, given three successive sampling instants $(t_0,
t_1, t_2)$ and their corresponding vectors $(Y_0, Y_1, Y_2)$, the
Fig. \ref{figure:headings1} can be interpreted as
follows.\begin{description}
\item \textit{Case a:} The sampling instants are chosen in such a
    way that the pair reachability /observability is preserved and
    also the pair controllability /constructibility .
\item \textit{Case b:} The sampling instants are chosen in such a
    way that the pair reachability /observability is not preserved
    but the pair controllability/ constructibility is preserved.
\item \textit{Case c:} The pair reachability /observability is
    neither preserved nor the pair controllability
    /constructibility. Remark that this situation will never occur
    for uniform sampling. Indeed, if the two first vectors are
    linearly dependent, then the third vector will always be
    linearly dependent.
\end{description}

From Case b, it can be seen that the pair controllability/constructibility
will always be guaranteed for any uniform sampling interval.
Nevertheless in the nonuniform case, this statement may not be true.

More intervals for the choice of the sampling instants in systems
of higher order can be found in \cite{Fuster1} and \cite{Troch}. According to the
previous example and references, the underlying idea of this note is
that, in spite of the restrictions of (17), there are intervals large
enough on the time axis where the sampling instants can be chosen
arbitrarily. Thus, this nearly free choice of the sampling instants can
be conveniently used in conjunction with other additional criterion
to optimize or improve different system aspects.

\section{Conclusions}
A joint characterization of reachability and observability for
linear SISO nonuniformly sampled discrete systems has been developed.
The classical characterization for uniformly sampled systems
appears as a simple particularization of the general criterion. The
nonuniform sampling offers a range of situations wider than the
uniform one in the study of these properties.

The note stresses the nonuniform sampling sequence, which is
believed to be an additional element for the analysis and handling of
discrete systems.

\end{document}